\documentstyle[12pt]{amsart}
\newtheorem{theo}{THEOREM}[section]
\newtheorem{lemma}[theo]{Lemma}
\newtheorem{cor}[theo]{Corollary}
\newtheorem{prop}[theo]{Proposition}
\theoremstyle{remark}
\newtheorem{rem}[theo]{Remark}

\newcommand{\brref}[1]{(\ref{#1})}
\newcommand{\restrict}[2]{{#1}_{\mid _{#2}}}
\newcommand{\Pin}[1]{{\Bbb P}^{#1}}
\newcommand{\lra}{\longrightarrow}
\title[Boundedness for Submanifolds of Quadrics]{Boundedness for Codimension Two Submanifolds of Quadrics}
\author[M.L. Fania ]{Maria Lucia Fania$^*$}
\address{ Maria Lucia Fania, Dipartimento di Matematica\\Universit\'{a} de L'Aquila\\
Via Vetoio Loc. Coppito\\67100 L'Aquila}
\email{fania@@univaq.it}
\author[G. Ottaviani] {Giorgio Ottaviani$^{*}$}
\address{Giorgio Ottaviani, Dipartimento di Matematica\\Universit\'{a} de L'Aquila\\
Via Vetoio Loc. Coppito\\67100 L'Aquila }
\email{ ottavian@@univaq.it}
\thanks{ $^*$ Both authors are members of GNASAGA of Italian CNR. 
Partially supported by MURST fund.}

\begin{document}
\maketitle
\vspace{5mm}
\begin{center}
{\it In ricordo di Fernando Serrano}
\end{center}
\vspace{5mm}

\begin{abstract}
It is proved that there are only finitely many families of codimension 
two subvarieties not of general type in  $Q_6$.
\end{abstract}

\section{introduction}
Ellingsrud and Peskine proved in \cite{ep} that 
smooth surfaces in $\Pin{4}$ not of general 
type have bounded degree. In \cite{boss1} this result 
has been estended to any non general type codimension two
submanifolds in $\Pin{n+2}$, of dimension $n \geq 2$.

In the same spirit Arrondo, Sols and De Cataldo proved in
\cite{as}, \cite{deC} the following result

\begin{theo}(Arrondo, Sols, De Cataldo) Let $X=X_n \subset 
Q_{n+2}$ be a smooth variety not of general type of dimension $n$ embedded in the
smooth quadric $Q_{n+2}$ of dimension $n+2$. Let $n \geq 2$, $n \neq 4$. 
Then deg(X) is bounded. 
\end{theo}

More precisely in \cite{as} it is proved the case n=2 while in 
\cite{deC} it is proved the case n=3 and it is observed that
the case $n \geq 5$ follows by an inequality on Chern Classes 
along the lines of \cite{sc}.

The aim of this paper is to drop the assumption $n \neq 4$ from the previous 
theorem. In fact we show the following

\begin{theo} Let $X=X_4 \subset 
Q_6$ be a smooth 4-dimensional variety not of general type. 
Then deg(X) is bounded. 
\end{theo}

As it is well known, see Prop. \brref{Hilbertcomp}, the theorem implies that there are only
finitely many families of codimension two subvarieties not of
 general type in $Q_6$. 

The paper is structured as follows.
In section 2. we fix our notation and give preliminary results
that will be needed later on in the paper. The sections 3. and 4.
are devoted to bounding the degree of a non general type
4-fold  $X \subset Q_6$.

In the last section we consider the problem of boundedness 
of non general type 4-folds in $\Pin{7}$.
Among the log-special type 4-fold in $\Pin{7}$ (i.e., the 
image of the adjunction  mapping has dimension less than 4)
there are still two hard cases to be considered, namely 
the quadric bundles over surfaces and the scrolls over threefolds.
We show that quadric bundles over surfaces have bounded 
degree with the only exception of those which lie in a 5-fold of degree 8.
 
The same technique can be applied to a manifold $X$ of 
dimension $n+1$ embedded in $\Pin{2n+1}$ which is a 
quadric bundle over a surface. One can prove that
there exists a function F(n) such that 
$deg(X) \leq F(n)$ or $X$ is contained in a variety of dimension 
$n+2$ and degree $[\frac{n(n+1)(4n^2+4n+1)}{6(4n-1)}]$,
where $[x]$ is the greatest integer less than or equal to x.

\section{Notations and Preliminaries}
\subsection{NOTATION}
\label{notation}
Throughout this article, unless otherwise specified, $X$ denotes  a smooth connected 
projective 4-fold defined over the complex field {\bf C}, 
which is contained in $Q_6$. 
Its structure sheaf is denoted  by ${\cal O}_X$.
For any coherent sheaf $\Im$ on $X$,  $h^i(\Im )$ is the complex dimension of
$H^i(S,\Im)$ and
$\chi=\chi({\cal O}_X)=\sum_i(-1)^ih^i({\cal O}_X).$
The following notation is used:\\
X, smooth 4-fold in  $Q_6$;\\
H class of hyperplane section of X, H =
$\restrict{{\cal O}_{\Pin{7}}(1)}{X}$;\\
K class of canonical bundle of X;\\
$X^3$ generic 3-fold section of X;\\
S generic surface section of X;\\
C generic curve section of X;\\
g genus of C;\\
$c_i$ = Chern classes of X;\\
N the normal bundle of X in $Q_6,  N_{X/Q_{6}}$.\\
Using the self intersection formula for the embedding of X $\subset Q_6$ 
\begin{equation}
\label{selfint}
c_2(N)=\frac{1}{2}dH^2
\end{equation}
we get the following formulae for $KH^3, K^2H^2, K^3H, K^4$ as 
function of d, g, $\chi({\cal O}_X)$, $\chi({\cal O}_{{X}^3})$,
$\chi(\cal{O}_S)$.

\begin{eqnarray}
\label{khcubo}
K\cdot H^3 = 2g-2-3d
\end{eqnarray}
\begin{eqnarray}
\label{kquadrohquadro}
K^2\cdot H^2 =
6\chi({\cal O}_S)-12g+12+\frac{13}{2}d+\frac{1}{4}d^2
\end{eqnarray}
\begin{eqnarray}
\label{kcuboh}
K^3\cdot H &=& -24\chi({\cal O}_{{X}^3})-48\chi({\cal O}_S)
+48g-48+3d-3d^2+\\
& & d(g-1)\nonumber\end{eqnarray}
\begin{eqnarray}
\label{kquarta}
 K^4 &=& 120\chi({\cal O}_X)+216\chi({\cal O}_{X^{3}})+
\chi({\cal O}_S)\frac{9d+472}{2}+ \\
& &\frac{5d^3+1098d^2-16d(45g+434)-6144(g-1)}{48}. \nonumber
\end{eqnarray}

Note that \brref{khcubo} follows from the adjuction 
formula.\\

To prove \brref{kquadrohquadro} we reason as follows. 
From the long exact sequence

\begin{equation}
\label{tangX}
0\lra T_X \lra T_{Q_{6}|X} \lra N \lra 0
\end{equation}
we get that
$$c_2(N) = 16H^2+6H\cdot K+K^2-c_2$$
Such equality along with the self intersection formula 
\brref{selfint} gives

\begin{equation}
\label{cdue}
c_2=(16-\frac{1}{2}d)H^2+6H\cdot K+K^2
\end{equation}
Hence by dotting \brref{cdue} with $H^2$ we have 
\begin{equation}
\label{cduehquadro}
c_2\cdot H^2=(16-\frac{1}{2}d)H^4+6H^3\cdot K+H^2\cdot K^2
\end{equation} 
On the other hand 
\begin{equation} 
\label{cduehduefrom}
c_2\cdot H^2=12\chi({\cal O}_S)-K^2\cdot H^2-(12g-12-11d).
\end{equation}
In fact using the following exact sequences

$$0\lra T_{X^{3}} \lra T_{X|{X^{3}}} \lra H_{X^{3}} \lra 0$$

$$0\lra T_{S}\lra T_{X^{3}|S} \lra H_S \lra 0$$
 we get that
\begin{equation}
\label{cdueh}
c_2\cdot H=c_2(X^3)-K_{X^{3}}\cdot H_{X^{3}}
\end{equation}
Now \brref{cduehduefrom} will follow if we 
dot \brref{cdueh} with H and if we use the following facts:

\begin{eqnarray}
\label{cdueXtre}
c_2(X^3)\cdot H_{X^{3}}&=&c_2(S)-K_S\cdot H_S=12\chi({\cal O}_S)-K_S^2-\\
& &(2g-2-d)\nonumber
\end{eqnarray}
and 
\begin{equation}
\label{ksdue}
K_S^2 = K^2\cdot H^2+8g-8-8d
\end{equation}

Combining \brref{cduehquadro} and \brref{cduehduefrom} 
we get  \brref{kquadrohquadro}.

In order to prove \brref{kcuboh} we do the following.
By dotting \brref{cdue} with H$\cdot$ K we get
\begin{equation} 
\label{cduehk}
c_2\cdot H\cdot K=(16-\frac{1}{2}d)K\cdot H^3+6K^2\cdot H^2+K^3\cdot H
\end{equation}
On the other hand if we dot \brref{cdueh} with K and we use
the fact that $c_2(X^3)\cdot K_{X^3}=-24\chi({\cal O}_{X^{3}})$ 
we have

\begin{equation} 
\label{cduehkfrom}
c_2\cdot H\cdot K=-24\chi({\cal O}_{X^{3}})-12\chi({\cal O}_S)+8g-8-6d
\end{equation}
 Now   \brref{kcuboh} is gotten by putting together \brref{cduehk} 
and \brref{cduehkfrom}.

We now prove  \brref{kquarta}. 
From the sequence \brref{tangX} we get that  
\begin{gather} 
\label{hcubo}
24H^3=c_3+c_2\cdot (6H+K)-\frac{1}{2}dK\cdot H^2\\
\label{hquarta}
22H^4=c_4+c_3\cdot (6H+K)-\frac{1}{2}dH^2\cdot c_2
\end{gather}
Thus 
\begin{eqnarray} 
\label{ctre}
c_3=(3d-72)H^3-(52-d)H^2\cdot K-12H\cdot K^2-K^3
\end{eqnarray}
\begin{eqnarray}
\label{cquattro}
c_4&=&454d-26d^2+\frac{1}{4}d^3+(384-12d)H^3\cdot K+\\
&  &\mbox{}(124-\frac{3}{2}d)H^2\cdot K^2+18H\cdot K^3+K^4\nonumber
\end{eqnarray}

By Riemann Roch theorem we know that 
\begin{equation} 
\label{RR}
-720\chi({\cal O}_X)=K^4-4K^2\cdot c_2-3c_2^2+K\cdot c_3+c_4
\end{equation}
Combining \brref{cdue}, \brref{RR}, \brref{cquattro}, \brref{ctre}, \brref{kcuboh}, 
\brref{kquadrohquadro} and \brref{khcubo}  we get \brref{kquarta}.

\begin{theo} (\cite{ccd}, Theorem 5.1)
\label{ccdbound} 
Let C be an irreducible reduced curve of arithmetic genus 
g and degree d in the  projective space $\Pin{n+1}$. Assume 
that C is not contained on any surface of degree  $<$ s, 
with $d > \frac{2s}{n-1}\prod_{i=1}^{n-1} \sqrt[n-i]{n!s}$. 
Then 
\begin{eqnarray*}
g-1 &\leq &\frac{d(d-1)}{2s}+\frac{d(s-2n+1)}{2(n-1)}+\frac{(d+s-1)(s-1)}{2s}+\\
& &\frac{(d-1)(n-2)(s+n-2)}{2s(n-1)}+\frac{(s-1)^2}{2(n-1)}
\end{eqnarray*}
\end{theo}

\begin{prop} (\cite{as}, Proposition 6.3)
\label{genusinQ3} 
Let C be a smooth curve of degree d and genus g in $Q_3$ 
that is not contained in any surface in  $Q_3$ of degree 
strictly less than 2k. Then 
$$g-1 \leq \frac{d^2}{2k}+\frac{1}{2}(k-4)d$$
\end{prop}

\begin{theo} (Castelnuovo bound \cite{har})
\label{harrisbound} 
Let V be an irreducible nondegenerate variety of 
dimension k and degree d in $\Pin{n}$. Put
\begin{xalignat}{2}
M=\left[\frac{d-1}{n-k}\right] &\quad and \quad \varepsilon = d-1-M(n-k)  \notag 
\end{xalignat} 
where $[x]$ is the greatest integer less than or equal to x. 
Then
\begin{eqnarray*}
p_g(V)= h^{0}(\tilde{V}, \Omega^{k}) \leq \binom{M}{k+1}(n-k)+\binom{M}{k}\varepsilon
\end{eqnarray*}
where $\widetilde{V}$ is a resolution of V (i.e., $\widetilde{V}$ is a smooth variety mapping 
holomorphically and birationally to V)
\end{theo}

%
%
\begin{prop} 
\label{Hilbertpoli}
 Let X be a smooth 4-fold in $Q_6$. Then
\begin{eqnarray*}
\chi({\cal O}_X(t)) &=&\frac{1}{24}dt^4+\frac{1}{12}(2-2g+3d)t^3+
\frac{1}{24}(12\chi({\cal O}_S)-12g+12+11d)t^2\\
&  &\mbox{}+\frac{1}{12}(12\chi({\cal O}_{X^{3}})+6\chi({\cal O}_S)-4g+4+3d)t+
\chi({\cal O}_X) 
\end{eqnarray*}

\end{prop}
\begin{pf} By the Riemann-Roch theorem we have 
$$\chi({\cal O}_X(t))= 
\frac{1}{24}H^4t^4+\frac{1}{12}(KH^3)t^3+
\frac{1}{24}(H^2\cdot K^2+c_2\cdot H^2)t^2-
\frac{1}{24}c_2\cdot H\cdot Kt+\chi({\cal O}_X),$$
where $c_i = c_i(T_X)$. We now use \brref{cduehduefrom} 
and  \brref{cduehkfrom} to get our claim.
\end{pf}

%
%
\begin{prop} 
\label{pgS}
 Let $S\subset Q_4$ be a surface of degree d contained 
in an irreducible threefold of degree $\sigma$, with $\sigma$ minimal.
Then
$$p_g(S)=h^2({\cal O}_S)\leq \frac{d^3}{24{\sigma}^2}+
\frac{d^2(\sigma-4)}{8\sigma}+
\frac{d(2{\sigma}^2-12\sigma+23)}{12}.$$
\end{prop}
\begin{pf} Let C be the generic curve section of S. By 
\cite{as}, Proposition 6.4  for $d>>0$ we have
$$g-1 \leq \frac{d^2}{4\sigma}+\frac{1}{2}(\sigma-3)d.$$
We let $$G(t) = \chi({\cal O}_{Q_{4}}(t))
=\left(\begin{array}{c}t+5\\5\end{array}\right)-
\left(\begin{array}{c}t+3\\5\end{array}\right)$$ and
$$\widetilde{F}(t) = G(t)-G(t-\sigma)-G(t-\frac{d}{2\sigma})+
G(t-\sigma-\frac{d}{2\sigma}).$$ We set

$$P(t) = dt+\left[-\frac{d^2}{4\sigma}+
\frac{1}{2}(3-\sigma)d \right]$$ and 

{$$\qquad{\text F(t)=}\left\{\aligned \widetilde{F}(t) \ \ \ \ \ if \ \ t \leq -1\\
\widetilde{F}(0)-1 \ \ \ if \ \ t=0\\0 \ \ \ if \ \ t\geq1 \endaligned \right.$$}

We have the following exact sequence 
 $$H^1({\cal O}_C(t))\lra H^2({\cal O}_S(t-1)) \lra H^2({\cal O}_S(t)) \lra 0$$
from which it follows that 

{$${-h^2({\cal O}_S(t))+
h^2({\cal O}_{S}(t-1))\leq h^1({\cal O}_C(t))=}
\left\{\aligned -\chi({\cal O}_C(t)) \leq -P(t)  \ \ \ \ \ if \ \ t \leq -1\\
-\chi({\cal O}_C)+1 \leq -P(0)+1  \ \ \ if \ \ t=0\endaligned\right\}$$}\\ = F(t-1)-F(t),  for t$\leq$ 0.\\

The same holds for t$\geq$ 1 since $h^1({\cal O}_C(t))=0$ 
and 
{$$\qquad{\text  F(t-1)-F(t)=}\left\{\aligned \geq 0 \ \ \ \ \ if \ \ t=1\\
0 \ \ \ if \ \ t \geq 2\endaligned \right.$$}
From this it  follows that F(t)-$h^2({\cal O}_S(t))$ is 
non increasing and since 
for t going to infinity it goes to zero it follows that 
$h^2({\cal O}_S(t)) \leq$ F(t) for all t. 
Thus evaluating it in t=0 we have
$$F(0)=\frac{d^3}{24{\sigma}^2}+
\frac{d^2(\sigma-4)}{8\sigma}+
\frac{d(2{\sigma}^2-12\sigma+23)}{12}$$ 
 and hence our claim. Note also that 
{$$\qquad{F(t)-F(t-1)=}
\left\{\aligned dt+[-\frac{d^2}{4\sigma}+
\frac{1}{2}(3-\sigma)d]  \ \ \ \ \ if \ \ t \leq -1\\
[-\frac{d^2}{4\sigma}+\frac{1}{2}(3-\sigma)d]-1 \ \ \ if \ \ t=0\endaligned\right.$$}\\

On passing note that if d is a multiple of 2$\sigma$ then
$\widetilde{F}(t)$ is the Hilbert polynomial of the complete intersection in $Q_4$ 
of hypersurfaces in  $\Pin{5}$ of degree $\sigma$ and 
$\frac{d}{2\sigma}$ while F(t) corresponds to 
$h^2({\cal O}_{V_{\sigma,\frac{d}{2\sigma}}}(t))$.
\end{pf}

%
%
\begin{prop}(\cite{deC}) 
\label{deC}
 Let $X^3 \subset V_{\sigma} \subset Q_5$. Then
$$-\chi({\cal O}_{X{^3}}) \geq  
\frac{1}{192{\sigma}^3}d^4+l.t. in \sqrt d.$$
\end{prop}
\begin{pf} For the proof see (\cite{deC}, Theorem 3.1.). 
It has to be noted that in there the coefficient 
 of $d^4$ should be $\frac{1}{192{\sigma}^3}$.
\end{pf}
%
%
\begin{prop} 
\label{Hilbertcomp}
For any fixed integer $d_0$ there are only finitely many 
irreducible components of the Hilbert scheme of  4-folds 
in $Q_6$ that contain 4-folds with $ d \leq d_0$.
\end{prop}
\begin{pf} By Harris' Castelnuovo bound, for $d \leq d_0$, 
there are finitely many possible values for g, $p_g(S), 
p_g(X^3), p_g(X)$ and hence for $\chi({\cal O}_S)$, 
$\chi({\cal O}_{X^{3}})$, $\chi({\cal O}_X)$
since $h^2({\cal O}_X) \leq p_g(S)$ and 
$h^1({\cal O}_S)=h^1({\cal O}_{X^{3}})=h^1({\cal O}_X)=0$. 
Thus there are only 
finitely many possibilities for the Hilbert polynomial
\begin{eqnarray*}
\chi({\cal O}_X(t))&=&\frac{1}{24}dt^4+\frac{1}{12}(2-2g+3d)t^3+
\frac{1}{24}(12\chi({\cal O}_S)-12g+12+11d)t^2\\
&  &\mbox{}+\frac{1}{12}(12\chi({\cal O}_{X^{3}})+6\chi({\cal O}_S)-4g+4+3d)t+
\chi({\cal O}_X) 
\end{eqnarray*}
\end{pf}

%
%
%
%
\section{4-Folds on a Hypersurface of Fixed Degree}

Let X be a 4-fold of degree d in $Q_6$ contained in an 
integral hypersurface $V_{\sigma} \in |{\cal O}_{Q_{6}}(\sigma)|.$

\begin{theo} 
\label{polin}
 Let $X \subset V_{\sigma} \subset Q_6$ be as above. There is a polynomial 
$P_{\sigma}(t)$ of degree 10 in $\sqrt d$ with positive leading coefficient, 
such that 
$$\chi({\cal O}_X) \geq P_{\sigma}(\sqrt d).$$
\end{theo}
\begin{pf} 
Look at the following three exact sequences

$$0 \lra {\cal O}_{\Pin{7}}(t-2) \lra 
{\cal O}_{\Pin{7}}(t) \lra 
{\cal O}_{Q_6}(t) \lra 0$$
$$0 \lra {\cal O}_{Q_6}(t-\sigma) \lra 
{\cal O}_{Q_6}(t) \lra 
{\cal O}_V(t) \lra 0$$
$$0 \lra {\cal I}_{X,V}(t) \lra 
{\cal O}_V(t) \lra 
{\cal O}_X(t) \lra 0$$

We use the first one to compute $\chi({\cal O}_{Q_6}(t))$ and 
the second one to compute

\begin{eqnarray*}
\label{chiOVt}
\chi({\cal O}_{V}(t))&=&\frac{{\sigma}}{60}t^5+
\frac{(6{\sigma}-{\sigma}^2)}{24}t^4+
\frac{\sigma({\sigma}^2-9\sigma+26)}{18}t^3-
\frac{\sigma({\sigma}^3-12{\sigma}^2+52\sigma-96)}{24}t^2\\
& &+\frac{\sigma(3{\sigma}^4-45{\sigma}^3+260{\sigma}^2-720\sigma+949)}{180}t\\
& &-\frac{\sigma({\sigma}^5-18{\sigma}^4+130{\sigma}^3-480{\sigma}^2+949\sigma-942)}{360}.
\end{eqnarray*}
Now use Prop.\brref{Hilbertpoli}, 
$\mu:=\mu_{\sigma}=\frac{1}{2}d^2+\sigma(\sigma-3)d-2\sigma(g-1)$
and the third exact sequence to compute

\begin{eqnarray*}
\label{chiideal}
\chi({\cal I}_{X,V}(t))&=&\frac{{\sigma}}{60}t^5+
\frac{1}{24}((6-{\sigma}){\sigma}-d)t^4+\\
& &[\frac{3d^2+6d\sigma(\sigma-6)-2(3\mu-2{\sigma}^2({\sigma}^2-9\sigma+26))}{72\sigma}]t^3-\\
& &[\frac{12\chi({\cal O}_S)\sigma-2d^2+d\sigma(19-4\sigma)+4\mu+{\sigma}^2({\sigma}^3-12{\sigma}^2+52\sigma-96)}{24\sigma}]t^2-\\
& & [\frac{180\chi({\cal O}_S)\sigma+360\chi({\cal O}_{X^3})\sigma-15d^2+30d\sigma(4-\sigma)+30\mu}{360\sigma}\\
& &-\frac{2{\sigma}^2(3{\sigma}^4-45{\sigma}^3+260{\sigma}^2-720\sigma+949)}{360\sigma}]t-\\
& &\frac{\sigma({\sigma}^5-18{\sigma}^4+130{\sigma}^3-480{\sigma}^2+949\sigma-942)}{360}-\chi({\cal O}_X)\\
&:=& Q(t)-\chi({\cal O}_X)
\end{eqnarray*}
Thus $\chi({\cal O}_X) = Q(t)-\chi({\cal I}_{X,V}(t))$.\\

Define 
\begin{eqnarray*}
t_1:= {\text {min}}\left \{ t\in {\Bbb N} | 
 \delta:=2\sigma t-d > 0,\frac{{\delta}^2}{2}-\mu-\delta\sigma(\sigma-3)>0\right\}
\end{eqnarray*}
Then
\begin{eqnarray*}
\frac{d}{2\sigma}\leq t_1\leq \frac{d}{2\sigma}+
\frac{\sqrt {2d}}{2}+\sigma
\end{eqnarray*}
By plugging $t_1$ we get that

\begin{eqnarray*}
Q(t_1)&\geq& \frac{1}{60\cdot 2^5{\sigma}^4}d^5-
\frac{1}{24\cdot 2^4{\sigma}^4}d^5-
\frac{1}{24\cdot 2^3{\sigma}^4}d^5-\\
& &\frac{1}{2^3{\sigma}^2}d^2 \chi({\cal O}_S)-
\frac{d}{2\sigma} \chi({\cal O}_{X^3})+ \textstyle{l.t. in \sqrt {d}}
\end{eqnarray*}
We now use Prop. \brref{deC} and Prop. \brref{pgS} in what above to get

\begin{eqnarray*}
Q(t_1)\geq \frac{d^5}{{\sigma}^4}(\frac{1}{60\cdot 2^5})+ \textstyle{l.t. in \sqrt {d}}
\end{eqnarray*}
and thus
\begin{eqnarray} 
 \chi({\cal O}_X) &\geq& Q(t_1)-\chi({\cal I}_{X,V}(t_1))\\
&\geq& -\chi({\cal I}_{X,V}(t_1))+
\frac{d^5}{{\sigma}^4}(\frac{1}{60\cdot 2^5})+ 
\textstyle{l.t. in \sqrt {d}}\nonumber
\end{eqnarray}
Moreover 
\begin{eqnarray*} 
 -\chi({\cal I}_{X,V}(t_1)) &\geq& -h^0({\cal I}_{X,V}(t_1))- 
h^2({\cal I}_{X,V}(t_1))-h^4({\cal I}_{X,V}(t_1))
\end{eqnarray*}
 
It will be enough to bound from above $h^2({\cal I}_{X,V}(t_1))$ since 
$h^0({\cal I}_{X,V}(t_1))$ and $h^4({\cal I}_{X,V}(t_1))$ 
have been bounded in (\cite {deC},  Lemma 3.3.).
In order to bound  $h^2({\cal I}_{X,V}(t_1))$ we consider 
the following exact sequences:

$$0 \lra {\cal O}_{Q_5}(-\sigma) \lra 
{\cal I}_{{X^3},{Q_5}} \lra 
{\cal I}_{{X^3},{V^4}} \lra 0$$
$$0 \lra {\cal O}_{\Pin{6}}(-2) \lra 
{\cal I}_{{X^3},{\Pin{6}}} \lra 
{\cal I}_{{X^3},{Q_5}} \lra  0$$ 

By \cite{bm}
\begin{xalignat}{2}
 H^1({\cal I}_{{X^3},{\Pin{6}}}(t)) = 0 &\quad for \quad t\geq4d-7 \notag 
\end{xalignat}
The latter along with the above sequences give  
\begin{xalignat}{2}
 H^1({\cal I}_{{X^3},{V^4}}(t)) = 0 &\quad for \quad t\geq4d-7 \notag 
\end{xalignat}
From 
$$0 \lra {\cal I}_{X,V}(k-1) \lra {\cal I}_{X,V}(k) \lra 
{\cal I}_{{X^3},{V^4}}(k) \lra  0$$
it follows that 
\begin{xalignat}{2}
 H^2({\cal I}_{X,V}(t)) = 0 &\quad for \quad t\geq 4d-8 \notag 
\end{xalignat}
Moreover by (\cite{deC}, Lemma 3.3) we have that
\begin{eqnarray*} 
 h^2({\cal I}_{X,V}(t_1)) \leq 
\sum_{k=t_1+1}^{4d-7} h^1({\cal I}_{{X^3},{V^4}}(k))
\leq (4d-7)A d^{\frac{7}{2}}+...
\leq 4Ad^{\frac{9}{2}}+...
\end{eqnarray*}

 Hence 
\begin{eqnarray*} 
 \chi({\cal O}_X) \geq Q(t_1)-\chi({\cal I}_{X,V}(t_1))
\geq Ad^5+ 
\mbox{...+ l.t. in d}
\end{eqnarray*}
which gives our claim.
\end{pf}

\begin{cor} 
\label{bound d}
  Let $X \subset V_{\sigma} \subset Q_6$ be as above. Assume that X is not of general type. 
Then there exists $d_0$ such that deg(X) $\leq d_0$.
\end{cor}
\begin{pf} Since X is not of general type we have 
$h^0(K_X(-1))$ = 0 from which it follows that 
$p_g(X) \leq p_g(X^3)$. This along with Harris bound give

\begin{eqnarray}
\label{chileq}
\chi({\cal O}_{X})&=&1+h^2({\cal O}_X)-h^3({\cal O}_X)+
p_g(X) \leq 1+h^2({\cal O}_X)\\
& &+p_g(X^3) \leq 1+p_g(S)+p_g(X^3)\leq  \frac{1}{216}d^4+ \mbox{l.t. in d}\nonumber 
\end{eqnarray}
On the other hand by \brref{polin} we obtain that 
\begin{eqnarray}
\label{chigeq}
\chi({\cal O}_X)\geq \frac{1}{1920{\sigma}^4}d^5+ &\mbox{l.t. in d} 
\end{eqnarray} 
The boundedness of d will now follow from \brref{chileq} 
and \brref{chigeq}. Hence our claim.
\end{pf}

%
%
\section{Boundedness}

\begin{prop} 
\label{chisless}
 Let X be a smooth 4-fold in $Q_6$. Denote 
$\chi({\cal O}_S)$, $\chi({\cal O}_{X^3})$,
$\chi({\cal O}_X)$ by s, x and v respectively. Then
\begin{itemize} 
\item[a)] $s \leq
\frac{2}{3}\frac{(g-1)^2}{d}+\frac{5}{3}(g-1)-\frac{1}{24}d^2+\frac{5}{12}d.$

\item[b)] $-24x(2g-2-3d) \leq
36s^2+3s(d^2-22d-16(g-1))+\frac{1}{16}(d^4-92d^3+4d^2(12g+193)+32d(1-g)(g+8)+768(g-1)^2).$
\end{itemize} 
If X is not of general type then 

\begin{itemize} 
\item[c)] $v \leq 2s-x$
\end{itemize}
\end{prop}
\begin{pf} By the generalized Hodge index theorem we know that 
\begin{gather}
\label{hodgeuno}
(K^2\cdot H^2)H^4 \leq (K\cdot H^3)^2\\
\label{hodgedue}
(K\cdot H^3)(K^3\cdot H) \leq (K^2\cdot H^2)^2.
\end{gather}

We observe that \brref{khcubo}, \brref{kquadrohquadro} and \brref{hodgeuno}
give a) while  \brref{khcubo}, \brref{kquadrohquadro}, 
\brref{kcuboh} and \brref{hodgedue} give b).
In order to prove c) we use the following exact sequence

\begin{equation}
\label{Ksequence}
0\lra K_{X}(-1)\lra K_X\lra K_{X^3}(-1)\lra 0
\end{equation}

Since X is not of general type we have $h^0(K_X(-1))$ = 0 
and thus $h^0(K_X) \leq h^0(K_{X^3}(-1)) \leq h^0(K_{X^3}).$
By the Lefschetz theorem it follows that $h^2({\cal O}_X) \leq h^2(\cal{O}_S)$. Moreover being
$h^1({\cal O}_X)$ = 0 it follows that $\chi({\cal O}_X)
\leq \chi({\cal O}_S)+p_g(X^3) =
2\chi({\cal O}_S)-\chi({\cal O}_{X^3}).$ Hence our claim.
\end{pf}
\begin{prop}
\label{chisgreater}
 Let X be a smooth 4-fold in $Q_6$. Then
\begin{itemize} \item[a)] $24\chi({\cal O}_S) \geq
d^2-2d-24(g-1).$
\item[b)] $240\chi({\cal O}_{X^3}) \leq
120\chi({\cal O}_X)+\frac{1}{2}(280-9d)\chi({\cal O}_S)-
\frac{1}{48}(d^3-6d^2+16d(12g-13)-960(g-1)).$
\end{itemize}  
\end{prop}
\begin{pf}
Since N(-1) is globally generated, the Segre classes satisfy:
\begin{xalignat}{2}
s_2(N(-1)) \cdot H^2 \geq 0, &\quad  s_4(N(-1)) \geq 0  \notag 
\end{xalignat} 
Recall that\\
$s_2 = c_1^2-c_2$\\ 
$s_4 = c_1^4+c_2^2-3c_1^2$.\\
Moreover\\
$c_1(N(-1)) = K+4H\\
c_2(N(-1)) = (\frac{1}{2}d-5)H^2 -
H\cdot K$.
Hence by \brref{khcubo}, \brref{kquadrohquadro}
\begin{eqnarray}
\label{stwo}
0&\leq&s_2(N(-1))\cdot H^2=K^2 \cdot H^2+16d+9K\cdot
H^3-(\frac{1}{2}d-5)H^4\nonumber \\
&= &6\chi({\cal O}_S)+6(g-1)+\frac{1}{2}d-\frac{1}{4}d^2
\end{eqnarray}
\begin{eqnarray}
\label{sfour}
0&\leq&s_4(N(-1))=(434-13d)K\cdot H^3+(136-\frac{3}{2}d)K^2\cdot H^2\\
& &+19K^3\cdot H+K^4+521d-29d^2+\frac{1}{4}d^3=-240\chi({\cal O}_{X^3})\nonumber \\ 
& &+120\chi({\cal O}_X)+\frac{280-9d}{2}\chi({\cal O}_S)-\frac{d^3-6d^2+16d(12g-13)}{48}\nonumber\\
& &-20(g-1)\nonumber
\end{eqnarray}
\end{pf}

\begin{theo}
\label{comphilb}
There are only finitely many irreducible components of the Hilbert scheme of smooth 4-folds in $Q_6$ that are 
not of general type.
\end{theo}
\begin{pf}
Let X be a smooth 4-fold in $Q_6$ that is not of general 
type. By Prop.\brref{Hilbertcomp} it is enough to bound d = degX. We will do so by 
considering separately the cases 2g-2-3d $\leq$ 0 and 2g-2-3d $>$ 0.
 
Assume that 2g-2-3d $\leq$ 0, i.e. g-1 $\leq \frac{3}{2}d$. 
Using Prop. \brref{chisless} along with a) in Prop. \brref{chisgreater} we get
$$0 \leq -\frac{1}{2}d^2+6d$$
Hence d is bounded in this case.

 Assume now that 2g-2-3d $>$ 0.
Using b) in Prop. \brref{chisgreater} along with c) and b) in 
Prop. \brref{chisless} we get 
\begin{eqnarray*}
 0 &\leq & \frac{540}{2g-2-3d}s^2 + \frac{117d^2-6d(3g+707)+80(g-1)}{2(2g-2-3d)}s +\frac{d^4}{2g-2-3d}\\
& &+\frac{-d^3(g+2078)+6d^2(229g+2842)+304d(1-g)(3g+23)}{24(2g-2-3d)}\\
& &+\frac{76(g-1)^2}{2g-2-3d}
\end{eqnarray*} 
Solving the above inequality with respect to s we see that either 
\begin{eqnarray}
\label{sgreater}
s \geq \frac{b+\sqrt {L}}{2160} \geq \frac{b}{2160}
\end{eqnarray} 
or 
\begin{eqnarray}
\label{sless}
s \leq \frac{b-\sqrt {L}}{2160}
\end{eqnarray}
where
$b = -117d^2+6d(3g+707)-80(g-1)$ and 

\begin{eqnarray*} 
L &=& 5049d^4-36d^3(107g+6793)+36d^2(9g^2-8978g+328809)+\\ 
& &960d(g-1)(339g+1915)-6560000(g-1)^2
\end{eqnarray*}
If \brref{sgreater} holds then combining it with 
Prop. \brref{chisless} we get 
\begin{equation}
\label{dquadro} 
0 \leq \frac{d^2}{80}- \frac{d(3g+557)}{360}+ \frac{2(g-1)^2}{3d)}- \frac{46}{27}(g-1).
\end{equation}
If \brref{sless} holds then such inequality along with a) 
in Prop. \brref{chisgreater} gives
\begin{eqnarray}
\label{dquarta} 
0&\leq& \frac{7}{864}d^4-\frac{d^3(5g+2203)}{6480}+ \frac{d^2(595-263g)}{3240}+\frac{7}{3}(g-1)^2\\
& &+\frac{d(1-g)(87g-5749)}{1620)}.\nonumber
\end{eqnarray}
Fix a positive integer k and let d $> 2k^2$. Assume that 
X does not lie on any hypersurface of $Q_6$ of degree 
strictly less than 2k. Then by Prop. \brref{genusinQ3} the genus of a general curve 
section of X satisfies 
\begin{equation}
\label{gless} 
g-1 \leq \frac{d^2}{2k} + \frac{1}{2}(k-4)d
\end{equation}
 
Rewriting \brref{dquadro} in the following way
$$0 \leq (g-1)\left[ \frac{2}{3d}(g-1)-\frac{d}{120}+\frac{46}{27} \right ]- \frac{7}{4}d + \frac{d^2}{80}$$  
and using \brref{gless} we get
\begin{eqnarray}
\label{glessk} 
(g-1) \leq \frac{3k}{2(k-40)}d + &\mbox {l.t. in d}.
\end{eqnarray} 
In the case \brref{dquarta} a  similar reasoning yields
\begin{eqnarray}
\label{glessd} 
(g-1) \leq \frac{21}{2}d + &\mbox{ l.t. in d}.
\end{eqnarray}
The following inequality, gotten by combining a) in  
Prop. \brref{chisless} 
and a) in Prop. \brref{chisgreater}  will be needed:
\begin{equation}
\label{eqcong} 
0 \leq -\frac{1}{12}d^2 + \frac{1}{2}d + (g-1)\left[\frac{2}{3d}(g-1)+\frac{8}{3}\right ]
\end{equation}
Plugging \brref{glessk} and \brref{gless} in 
\brref{eqcong} gives 
\begin{eqnarray}
\label{i} 
0 \leq d^2(\frac{1}{2(k-40)}-\frac{1}{12}) + &\mbox{l.t. in d}.
\end{eqnarray}
Similarly, plugging \brref{glessd} and \brref{gless} in 
\brref{eqcong} gives
\begin{eqnarray}
\label{ii} 
0 \leq d^2(\frac{7}{2k}-\frac{1}{12}) + &\mbox{l.t. in d}.
\end{eqnarray}
The coefficient of $d^2$, both in \brref{i} and \brref{ii} 
is negative for k=47. Hence d is bounded from above if X is not 
in a hypersurface of degree strictly less than 2 $\cdot$ 47.
If X is not of general type and is contained in a 
hypersurface of degree less than or equal to 2$\cdot$ 47 
then by Corollary \brref{bound d} there exists $d_0$ such that 
deg(X) $\leq d_0$. Hence the theorem is proved.
\end{pf}

%
%

\section{Quadric bundles over surfaces in $\Pin{7}$}

Throughout this section $X$ will denote  a smooth 4-fold 
of degree d in $\Pin{7}$ which is a quadric bundle over 
a surface. We will show that either its degree d is bounded 
or $X$ is contained in a $5$-fold of degree $8$.

For 4-folds in $\Pin{7}$ by the selfintersection formula we have:
$$c_3(N_{X|\Pin{7}}) = dH^3$$

\subsection{}
\label{ciXinP7}  
From the exact sequence
  
\begin{eqnarray*}
\label{tangXinP7}
0\lra T_X \lra T_{\Pin{7}|X} \lra N_{X|\Pin{7}} \lra 0
\end{eqnarray*}
we get that

\begin{eqnarray*}
c_3(X) &=& (56-d)H^3-28H^2\cdot c_1(X)+
           8H\cdot(c_1^2(X)-c_2(X))-c_1^3(X)+\\
& &2c_1(X)c_2(X)\end{eqnarray*}

\begin{eqnarray*}
c_4(X) &=& 70d-c_1(X)dH^3-c_2(X)(28H^2-8c_1(X)H+c_1^2(X)-c_2(X))\\
& &-c_3(X)(8H-c_1(X))\end{eqnarray*}

\subsection{ Definition} A 4-fold $X$ is called a geometric 
quadric bundle if there exists a morphism $p:X \lra B$ 
onto a normal surface B such that every fibre
$p^{-1}(b)$ is isomorphic to a quadric.  A 4-fold X is a 
quadric bundle in the adjuction theoretic sense 
if there exists a morphism $p:X \lra B$ onto a normal surface B 
and an ample Cartier divisor L on B such $p^{*}L=K+2H$.\\ 

The following proposition relates the two notions.
\begin{prop}
\label{...}
Let X  be a quadric bundle in the 
adjuction theoretic sense. Then X is a geometric quadric bundle.
Moreover the base B is smooth.
\end{prop}
\begin{pf}
By (\cite {beso1}, Theorem 2.3) we know that p is 
equidimensional, being dim X=4. Moreover by  
(\cite {bes}, Theorem 8.2) the base B is smooth.
\end{pf}

We fix our notation which follows closely the one in 
\cite{boss2}.

\subsection{Notations}
\label{notazione quadric fib}  Let $p:X \lra B$  be a 
geometric quadric bundle in  $\Pin{7}$. 
We have a natural morphism $f:B \lra Gr(\Pin{3},\Pin{7})$.
   
Let S be a generic surface section of X. Then 
$p:S \lra B$ is finite 2:1. Let 2R $\subset B$ be the ramification divisor of $p:S \lra B$.

 We set  $p_{*}{{\cal O}_{X}}(1)$:= E, a rank 4 vector 
bundle over B. We have E = $f^{*}(U^{\vee})$,
where U is the universal bundle of  $Gr(\Pin{3},\Pin{7})$, 
in particular det E = $f^{*}(\stackrel{2}{\wedge}U^{\vee})$ 
is ample. Note that W:= P(E) is a $\Pin{3}$-bundle in the natural 
incidence variety $\Pin{7}\times Gr(\Pin{3},\Pin{7})$ whose 
projection $\pi$ into $\Pin{7}$ is the hypersurface V given
 by the union of all the 3-planes containing the quadrics 
of X.

Moreover $\pi^{-1}(X) = \widetilde{X}$ is smooth and isomorphic to X.
We denote the natural projection of W onto B also 
by p and by H the divisor on W corresponding to 
${\cal O}_{W}(1)$.
Hence  $\widetilde{X} = 2H-p^{*}L$ for some divisor L on B.
   
The divisor D $\subset B$ corresponding to points whose 
fibres are singular quadrics, is called the discriminant 
divisor. Moreover D = $c_1(E)-c_1(L\otimes E^{\vee}) 
= 2c_1(E)-4L.$ In fact $\widetilde{X}$ determines a section of
$S^{2}E\otimes L^{\vee}$, hence a morphism 
$\phi: L\otimes E^{\vee}\lra E$. D is given by the equation
 det$\phi$ = 0. Thus our claim.\\

In order to bound d we need several preliminaries computations.

\begin{prop}
\label{ciW}
\begin{eqnarray*}
c_1(W)&=&4H-p^{*}c_1(E)+p^{*}c_1(B)\\
c_2(W)&=&6H^2+H\cdot[4p^{*}c_1(B)-3p^{*}c_1(E)]+
p^{*}c_2(B)-p^{*}c_1(E)\cdot p^{*}c_1(B)\\
& &+p^{*}c_2(E)\\
c_3(W)&=&4H^3+H^2\cdot[6p^{*}c_1(B)-3p^{*}c_1(E)]+
H\cdot[2p^{*}c_2(E)+4p^{*}c_2(B)\\
& &-3p^{*}c_1(E)\cdot p^{*}c_1(B)] \\
c_4(W)&=&4H^3\cdot p^{*}c_1(B)+H^2\cdot[6p^{*}c_2(B)-3p^{*}c_1(E)\cdot 
p^{*}c_1(B)]
\end{eqnarray*}
$H^4-H^3\cdot p^{*}c_1(E)+H^2\cdot p^{*}c_2(E)=0$
\end{prop}

\begin{pf} Consider the sequence
\begin{eqnarray*}
\label{tangW}
0\lra {\cal O}_{W} \lra p^{*}E^{\vee} \otimes {\cal O}_{W}(1) \lra T_W \lra p^{*}T_B \lra 0
\end{eqnarray*}
The Chern polynomial of 
$p^{*}E^{\vee} \otimes {\cal O}_{W}(1)$ 
is 

$1+c_1(p^{*}E^{\vee} \otimes {\cal O}_{W}(1))t+
c_2(p^{*}E^{\vee} \otimes {\cal O}_{W}(1))t^2+
c_3(p^{*}E^{\vee} \otimes {\cal O}_{W}(1))t^3 = 
1+[4H-p^{*}c_1(E)]t+[6H^2-3p^{*}c_1(E)\cdot H+
p^{*}c_2(E)]t^2+[4H^3-3p^{*}c_1(E)\cdot H^2+2
p^{*}c_2(E)\cdot H]t^3.$

 On the other hand
$ch(T_W) = ch(p^{*}E^{\vee} \otimes {\cal O}_{W}(1))
\cdot ch(p^{*}T_B)$
hence we get that
$1+c_1(W)t+c_2(W)t^2+c_3(W)t^3+c_4(W)t^4+c_5(W)t^5=\{1+
[4H-p^{*}c_1(E)]t+[6H^2-3p^{*}c_1(E)\cdot H+
p^{*}c_2(E)]t^2+[4H^3-3p^{*}c_1(E)\cdot H^2+2
p^{*}c_2(E)\cdot H]t^3\}\cdot \{1+p^{*}c_1(B)t+p^{*}c_2(B)t^2\}$.
Expanding the right hand side we get the first four equations.
The last one is the Wu-Chern equation on W = P(E), that is,
$c_4(p^{*}E^{\vee} \otimes {\cal O}_{W}(1))=0.$
\end{pf}  

\begin{lemma} 
\label{c1E}
$c_1(E) = 2R-\frac{D}{2}, \quad L = R-\frac{D}{2}$.
\end{lemma}
\begin{pf} We have $K_S=p^{*}(K_B+R)$, 
hence by the adjunction formula 
$K_X=-2H+p^{*}R+p^{*}K_B$.
From Prop. \brref{ciW} $K_W=-4H+p^{*}c_1(E)-p^{*}c_1(B)$. Putting 
this together with the adjunction formula
$K_X=\restrict{K_W}{X}+2H-p^{*}L$ gives
$-p^{*}L+p^{*}c_1(E) = p^{*}R$, that is $c_1(E)=L+R$.
Substituting this in $c_1(E)=\frac{D+4L}{2}$ we get 
$L = R-\frac{D}{2}$ and hence $c_1(E) = 2R-\frac{D}{2}$.
\end{pf}

\begin{prop} 
\label{ciX}
\begin{eqnarray*}
c_1(X) &=& 2H-p^{*}K_B-p^{*}R\\
c_2(X) &=& 2H^2+H\cdot[-p^{*}c_1(E)+2p^{*}c_1(B)]+p^{*}{R^2}-
p^{*}c_1(E)\cdot p^{*}R\\
& &-p^{*}c_1(B)\cdot p^{*}R+ p^{*}c_2(B)+p^{*}c_2(E)\\
c_3(X) &=& H^2\cdot [2p^{*}c_1(B)+p^{*}c_1(E)-
2p^{*}R]+ H\cdot[2p^{*}c_2(B)-\\
& &p^{*}c_1(B)\cdot p^{*}c_1(E)-2p^{*}{R^2}+3p^{*}R\cdot p^{*}c_1(E)-p^{*}c_1^2(E)]\\
c_4(X) &=& H^3\cdot[-2p^{*}c_1(E)+4p^{*}R]+H^2\cdot[2p^{*}c_2(B)+p^{*}c_1(B)\cdot p^{*}c_1(E)\\
& &+6p^{*}{R^2}+3p^{*}c_1^2(E)-9p^{*}R\cdot p^{*}c_1(E)-
2p^{*}R\cdot p^{*}c_1(B)]
\end{eqnarray*}
\end{prop}

\begin{pf} The following sequence
\begin{eqnarray*}
0\lra T_X \lra T_{W|X} \lra {\cal O}(2H+R-c_1(E)) \lra 0
\end{eqnarray*}
along with Prop. \brref{ciW} gives the proof.
\end{pf}  

\begin{lemma} Let Z, $Z'$ be arbitrary divisors on B. Then
\begin{itemize} 
\item [i)] $H^2\cdot p^{*}Z\cdot p^{*}Z' = 2Z\cdot Z'$

\item[ii)] $H^3\cdot p^{*}Z = (c_1(E)+R)\cdot Z 
= (3R-\frac{D}{2})\cdot Z$
\end{itemize}
\label{intersecform}
\end{lemma}
\begin{pf} i) follows from the fact that the fibres of X over B 
are quadrics. 
As for ii) note that $H^3\cdot p^{*}Z$ is equal to 
the intersection product in W
$$ H^3\cdot p^{*}Z\cdot (2H+p^{*}R-p^{*}c_1(E)).$$
Intersecting the Wu-Chern equation with 
$p^{*}Z$ we get that
$$ H^4\cdot p^{*}Z-H^3\cdot p^{*}Z\cdot p^{*}c_1(E) 
=0.$$ Hence

$H^3\cdot p^{*}Z\cdot (2H+p^{*}R-p^{*}c_1(E))=
2H^4\cdot p^{*}Z+ H^3\cdot p^{*}Z\cdot p^{*}R
- H^3\cdot p^{*}Z\cdot p^{*}c_1(E)
=2H^3\cdot p^{*}Z\cdot p^{*}c_1(E)+H^3\cdot p^{*}R\cdot p^{*}Z- H^3\cdot p^{*}Z\cdot p^{*}c_1(E)
=(c_1(E)+R)\cdot Z = (3R-\frac{D}{2})\cdot Z$
\end{pf}

\begin{prop} 
\label{deg,c2E}
The surface f(B) in $Gr(\Pin{3},\Pin{7})$ has bidegree 
$(\delta, c_2(E))$ where
\begin{xalignat}{2}
 \delta = degV = \frac{d-R\cdot c_1(E)+c_1^2(E)}{2}, &\quad   
c_2(E)= \frac{-d+R\cdot c_1(E)+c_1^2(E)}{2}\notag 
\end{xalignat}
\end{prop}
\begin{pf} We intersect the Wu-Chern equation in Prop.
\brref{ciW} with H and we get 
 $$\delta -H^4\cdot c_1(E)+H^3\cdot c_2(E) = 0.$$ Now cut 
the equation $X = 2H+p^{*}R-p^{*}c_1(E)$ with $H^4$ and
 we obtain 
$$d = 2\delta+H^3\cdot p^{*}R\cdot p^{*}c_1(E)-
H^3\cdot p^{*}c_1^2(E).$$
From these two equalities we get our claim.
\end{pf}  

\begin{prop} 
\label{ciXconsostituzioni}
\begin{eqnarray*}
c_1(X) &=& 2H+p^{*}c_1(B)-p^{*}R\\
c_2(X) &=& 2H^2+H\cdot[2p^{*}c_1(B)-2p^{*}R+\frac{1}{2}p^{*}D]
-p^{*}{R^2}+p^{*}c_2(E)-\\
& & p^{*}c_1(B)\cdot p^{*}R+p^{*}c_2(B)+\frac{1}{2}p^{*}D\cdot p^{*}R\\
c_3(X) &=& H^2\cdot[-\frac{1}{2}p^{*}D+
2p^{*}c_1(B)]+ H\cdot[\frac{1}{2}p^{*}c_1(B)\cdot p^{*}D-\\
& &2p^{*}c_1(B)\cdot p^{*}R+2p^{*}c_2(B)-\frac{1}{4}p^{*}{D^2}
+\frac{1}{2}p^{*}D\cdot p^{*}R]\\
c_4(X) &=& H^3\cdot p^{*}D+H^2\cdot[-\frac{1}{2}p^{*}c_1(B)\cdot p^{*}D
+2p^{*}c_2(B)+\frac{3}{4}p^{*}{D^2}-\frac{3}{2}p^{*}D\cdot p^{*}R]
\end{eqnarray*}
\end{prop}
\begin{pf} Using \brref{c1E} - \brref{deg,c2E}
 and easy computations give the formulas for 
$c_2(X), c_3(X), c_4(X)$.
\end{pf}

\begin{prop} 
\label{system}
Set $$P(d)=\frac{1}{9d^3-50d^2-10949d+169120}, \quad x = K_B^2, \quad y = D\cdot R.$$ 
The following hold:\\
\begin{eqnarray*}
v=R^2&=&-\frac{1}{2}P(d)[-192864x+2842{d}^{3}+332024d-70224y-
\\
& &53900{d}^{2}-15y{d}^{2}-49{d}^{4}+4974yd-3y{d}^{3}+9016dx]\\
z=D^2&=&P(d)[-882{d}^{4}+54y{d}^{3}+54684{d}^{3}-870y{d}^{2}-1100736{d}^{2}+\\
& &190512dx-39792yd+7112448d-3035648x+709632y]\\
f=c_2(B)&=&\frac {1}{16}P(d)[9{d}^{5}-328{d}^{4}+
144{d}^{3}x+3y{d}^{3}+3036{d}^{3}-173y{d}^{2}-\\
& &37416{d}^{2}+712{d}^{2}x+2608yd-102048dx+728896d+\\
& &675584x-8576y]\\
u=K_B\cdot R&=&-\frac{1}{8}P(d)[-175392dx+1722{d}^{4}+696696{d}^{2}
-22686yd-\\
& &52332{d}^{3}-11y{d}^{3}-21{d}^{5}+879y{d}^{2}+3864{d}^{2}x+1983744x-\\
& &3415104d+190784y]\\
t=K_B\cdot D&=&-\frac{1}{4}P(d)[-63{d}^{5}+5292{d}^{4}
-164556{d}^{3}-33y{d}^{3}+
2237760{d}^{2}+\\
& &13608{d}^{2}x+2703y{d}^{2}-11176704d-71880yd-516208dx+\\
& &624384y+4770304x]
\end{eqnarray*}
\end{prop}
\begin{pf} We get a linear system of five equations in 
the unknowns f, t, u, v, z with coefficients 
rational functions of d. In fact from \brref{ciXinP7} we 
have
$$c_3(X) = (56-d)H^3-28H^2\cdot c_1(X)+8H\cdot(c_1^2(X)
-c_2(X))-c_1^3(X)+2c_1(X)c_2(X).$$ Substituting the values of 
$c_1(X), c_2(X), c_3(X)$ of Proposition \brref{ciXconsostituzioni}
we get

$(d-16){H}^{3}+\left(-12p^*R+\frac{3p^*D}{2}+
14p^*c_1(B)\right)\cdot {H}^{2}+
(-6{p^*c_1(B)}^{2}-10{p^*R}^{2}+
6p^*c_2(B)+4p^*c_2(E)-\frac {{p^*D}^{2}}{4}+6p^*c_1(B)
\cdot p^*R-\frac{p^*c_1(B) \cdot p^*D}{2}+\frac{7p^*R\cdot p^*D}{2})\cdot H =0$

Cutting respectively with H, 
$p^{*}R$, $p^{*}K_B$, $p^{*}D$ we get four equations.
 For instance if we cut with 
$p^{*}D$ we obtain:
$$(d-16)H^3\cdot p^{*}D
+H^2\cdot \left (-12p^*R\cdot p^*(D)+\frac{1}{2}3p^*D\cdot p^*D
+14p^*c_1(B)\cdot p^*D\right )=0$$

Using now lemma \brref{c1E} and 
lemma \brref{intersecform} we get
$$(d-16)(3R-\frac{D}{2})\cdot D-24 R\cdot D+3D\cdot D+
28c_1(B)\cdot D=0$$

and simplifying

$$-72R\cdot D+11{D}^{2}+3dR\cdot D-\frac{{dD}^{2}}{2}
-28K_B\cdot D=0$$

The fifth equation is gotten by substituting the values 
of Proposition \brref{ciXconsostituzioni}
in the second formula of \brref{ciXinP7}:
\begin{eqnarray*}
c_4(X) &=& 70d-c_1(X)dH^3-c_2(X)(28H^2-8c_1(X)H+c_1^2(X)-c_2(X))\\
& &-c_3(X)(8H-c_1(X))\end{eqnarray*}
Solving such system we get the claim.
\end{pf}

\begin{prop} 
\label{c2E,g}
\begin{eqnarray*}
c_2(E)&=&\frac{1}{4}P(d)[\left(41160d-360640\right )x+\left (-95{d}^{2}+5005d-69440\right)y\\
& &-165{d}^{4}+10390{d}^{3}-205070{d}^{2}+1225840d]\\
g-1&=&\frac{1}{8}P(d)[\left(1008{d}^{2}-49112d+566720\right)x
+(223{d}^{2}-9857d+\\
& &109120)y+393{d}^{4}-21032{d}^{3}+353440{d}^{2}-1781360d]
\end{eqnarray*}
\end{prop}
\begin{pf} Using Prop. \brref{deg,c2E} and Lemma \brref{c1E} we get
$$c_2(E)=\frac{1}{2}\left(6R^2-\frac{5D\cdot R}{2}+\frac{D^2}{4}-d\right)$$
Moreover from the adjuction formula and Lemma \brref{intersecform}
ii) we obtain
\begin{eqnarray*}
g-1&=&\frac{1}{2}d+\frac{1}{2}H^3(p^{*}K_B+p^{*}R)=
\frac{1}{2}d+\frac{1}{2}(p^{*}K_B+p^{*}R)(3R-\frac{D}{2})\\
&=&\frac{1}{2}d+\frac{3}{2}K_B\cdot R-\frac{1}{4}K_B\cdot D+\frac{3}{2}R^2-\frac{1}{4}D\cdot R.
\end{eqnarray*}
Substituting the values of Prop. \brref{system} and simplifying 
we get the assertions.
\end{pf}

We need a Roth type result.

\begin{prop} Let X be a codimension $3$ integral subvariety
of $\Pin{n}$ of degree $d$. If the generic section 
$C=X\cap \Pin{4}$
with a linear $\Pin{4}$ 
is contained in a surface $S_{\sigma}\subset\Pin{4}$ of 
degree $\sigma$ with
$\sigma^2 \leq d$ then X itself is contained in a codimension
$2$ subvariety $V_{\sigma}\subset\Pin{n}$ of degree $\sigma$.
\label{roth}
\end{prop}
\begin{pf}
We first check that the generic section $S=X\cap \Pin{5}$
with a linear $\Pin{5}$ is contained in a 3-fold
$Y_{\sigma}\subset\Pin{5}$ of degree $\sigma$.
  By the assumptions, the surface $S_{\sigma}$ is unique.
  On the contrary, suppose there are two such surfaces
 $S'_{\sigma}$ and $S''_{\sigma}$.  
There exists a linear projection $\pi$ from
  a point $p\in\Pin{4}$ on a hyperplane $\Pin{3}$ such that
$\pi (S'_{\sigma})$ and $\pi (S''_{\sigma})$ are two
irreducible
distinct surfaces of degree $\sigma$ containing $\pi (C)$
against the Bezout theorem. We get a rational map
from  $(\Pin{5}) ^{\vee}$ in the Hilbert scheme
of degree $\sigma$ surfaces of $\Pin{4}$. It follows that
the closure of all surfaces $S_{\sigma}$ is the 3-fold
$Y_{\sigma}$ we looked for.

By iterating this process we get the thesis.
\end{pf} 

\begin{prop} Let X be a quadric bundle in $\Pin{7}$. Then 
$d\leq 2963$ or $X$ is contained in a $5$-fold of degree $8$.
\label{bound}
\end{prop}
\begin{pf} We consider the possible values of x and y compatible with
the following three inequalities:

$D\cdot R \geq 0$,  \quad\mbox{(Lemma 4.15 in \cite{boss2})}\\
$c_2(E) \geq 0$,  \quad\mbox{(Proposition \brref{deg,c2E})}\\
$c_1(E)\cdot D \geq 0$,  \quad\mbox{($c_1$(E) is ample, see 
Notations \brref{notazione quadric fib})}.\\

Using Lemma \brref{intersecform}, Prop. \brref{deg,c2E}, Prop. \brref{system},
Prop. \brref{c2E,g} and easy computations we have\\
$y \geq 0\\
\left (41160d-360640\right )x
+\left (-95{d}^{2}+5005d-69440\right 
)y-165{d}^{4}+10390{d}^{3}-205070{d}^{2}\\
+1225840d \geq 0\\
-\left (190512d-3035648\right )x
- \left (18{d}^{3}-670{d}^{2}+4004d
+33152\right )y+882{d}^{4}-54684{d}^{3}\\
+1100736{d}^{2}-7112448d
\geq 0$\\

These inequalities bound the inside of a triangle 
whose vertices are: 

$$A_d=        \left (\frac {d\left (
33{d}^{3}-2078{d}^{2}+41014d-245168\right )}{8232d-72128},0\right )$$
          
$$B_d=\left (\frac {9d\left ({d}^{3}-62{d}^{2}+1248d-8064
\right )}{1944d-30976},0\right )$$
                              
$$C_d=\left (\frac {\left (33{d}^{3}-2192{d}^{2}+46900d-316064  
\right )d}{8232d-131712},\frac {2d\left (69d-872\right )
}{21d-336}\right )$$

The minimum and the maximum of g-1 
considered as a function in x and y with d fixed (see Prop. \brref{c2E,g}) have to be assumed 
in one of the vertices. Substituting the coordinates of  
$A_d, B_d, C_d$ in the expression of g-1 we get 
respectively

$g-1=\frac{(33{d}^{2}-293d-552)d}{588d-5152}$,
$\frac{(63{d}^{2}-1320d+5192)d}{972d-15488}$,
$\frac{(33{d}^{2}-407d-1008)d}{588d-9408}.$
 
Hence we have in our case
\begin{eqnarray}
\label{gcompreso}
{\frac {\left (33{d}^{2}-293d-552\right )d}{588d-5152}}
\leq g-1 \leq{
\frac {\left (63{d}^{2}-1320d+5192\right )d}{972d-15488}}
\end{eqnarray}
We now distinguish two cases.
Suppose first that the curve section C is not contained in 
any surface of degree 8. Then from Theorem \brref{ccdbound} and
\brref{gcompreso} we have

$${\frac {\left (33{d}^{2}-293d-552\right )d}{588d-5152}}
\leq  
\frac{d^2}{18}+\frac{5}{3}d+\frac{347}{18}$$
that is $3d^3-8881d^2-29706d+893872 \leq 0$, which 
implies $d \leq 2963$.

If otherwise then C is contained in a octic, then from
Proposition \brref{roth} it follows that also $X$ is contained in a octic. 
\end{pf}

\begin{rem}
We like to remark that a similar reasoning gives an 
analogous result for manifolds $X$ of dimension $n+1$ and 
degree $d$ embedded in $\Pin{2n+1}$ which is a quadric bundle
over a surface. More precisely we have the following 
\end{rem}
\begin{prop} Let $X$ be a manifold of dimension $n+1$ and 
degree d embedded in $\Pin{2n+1}$ which is a quadric bundle
over a surface. Then there exists a function F(n) such that
$d \leq F(n)$ or $X$ is contained in a variety of dimension 
$n+2$ and degree $[\frac{n(n+1)(4n^2+4n+1)}{6(4n-1)}]$,
where $[x]$ is the greatest integer less than or equal 
to x.
\label{Bound}
\end{prop}
\begin{pf} 
Similarly as in proposition \brref{bound} we get that 
$$Ad^2+O(d) \leq g-1 \leq 
Bd^2+O(d)$$
where $A=\frac{3(4n-1)}{n(n+1)(4n^2+4n+1)}$ and 
$B=\frac{3(2n+1)}{(n+1)(2n^3+2n^2+2n+3)}$.
We now use  Theorem \brref{ccdbound} and Proposition \brref{roth} to 
see that there exists a function F(n) such that
$d \leq F(n)$ or $X$ is contained in a variety of dimension 
$n+2$ and degree $[\frac{n(n+1)(4n^2+4n+1)}{6(4n-1)}]$,
where $[x]$ is the greatest integer less than or equal 
to x. 
\end{pf}


\end{document}